\documentclass[a4paper,UKenglish,cleveref, autoref, thm-restate]{lipics-v2021}

\hideLIPIcs  

\title{Finding \texorpdfstring{$d$-Cuts}{d-Cut} in Claw-free Graphs} 

\author{Jungho Ahn}{Department of Computer Science, Durham University, Durham, UK}{jungho.ahn@durham.ac.uk}{https://orcid.org/0000-0003-0511-1976}{supported by Leverhulme Trust Research Project Grant RPG-2024-182.}
\author{Tala Eagling-Vose}{Department of Computer Science, Durham University, Durham, UK}{tala.j.eagling-vose@durham.ac.uk}{https://orcid.org/0009-0008-0346-7032}{}
\author{Felicia Lucke}{Department of Computer Science, Durham University, Durham, UK}{felicia.c.lucke@durham.ac.uk}{https://orcid.org/0000-0002-9860-2928}{supported by EPSRC Grant EP/X01357X/1.}
\author{Dani\"el Paulusma}{Department of Computer Science, Durham University, Durham, UK}{daniel.paulusma@durham.ac.uk}{https://orcid.org/0000-0001-5945-9287}{supported by grants RPG-2024-182 and EP/X01357X/1.}
\author{Siani Smith}{Department of Computer Science, Loughborough University, Loughborough, UK}{s.smith16@lboro.ac.uk}{https://orcid.org/0000-0003-0797-0512}{}

\authorrunning{J. Ahn, T. Eagling-Vose, F. Lucke, D. Paulusma and S. Smith}

\Copyright{Jungho Ahn, Tala Eagling-Vose, Felicia Lucke, Dani\"el Paulusma and Siani Smith}

\ccsdesc[500]{Mathematics of computing~Graph theory}
\ccsdesc[500]{Theory of computation~Graph algorithms analysis}
\ccsdesc[500]{Theory of computation~Problems, reductions and completeness}

\keywords{matching cut, \texorpdfstring{$d$-cut}{d-cut}, claw-free, maximum degree.}

\category{} 

\relatedversion{} 

\nolinenumbers 

\EventEditors{John Q. Open and Joan R. Access}
\EventNoEds{2}
\EventLongTitle{42nd Conference on Very Important Topics (CVIT 2016)}
\EventShortTitle{CVIT 2016}
\EventAcronym{CVIT}
\EventYear{2016}
\EventDate{December 24--27, 2016}
\EventLocation{Little Whinging, United Kingdom}
\EventLogo{}
\SeriesVolume{42}
\ArticleNo{23}

\newtheorem{open}{Open Problem}

\newcommand{\NP}{{\sf NP}}
\newcommand{\ssi}{\subseteq_i}
\newcommand{\si}{\supseteq_i}

\usepackage{tikz}
\usetikzlibrary{decorations.pathreplacing, calligraphy, math}

\definecolor{nicered}{RGB}{204,0,0}
\definecolor{lightblue}{RGB}{153,204,255}
\tikzstyle{vertex}=[thin,circle,inner sep=0.cm, minimum size=1.7mm, fill=black, draw=black]
 \tikzstyle{svertex}=[thin,circle,inner sep=0.cm, minimum size=1.3mm, fill=black, draw=black]
 \tikzstyle{bvertex}=[thin,circle,inner sep=0.cm, minimum size=1.7mm, fill=lightblue, draw=lightblue]
 \tikzstyle{rvertex}=[thin,circle,inner sep=0.cm, minimum size=1.7mm, fill=nicered,draw=nicered]
 \tikzstyle{evertex}=[thin,circle,inner sep=0.cm, minimum size=1.7mm, fill=none,draw=black]
 \tikzstyle{edge}=[thick, draw = gray]
 \tikzstyle{tedge}=[ultra thick, draw = black]
 \tikzstyle{tredge}=[ultra thick, draw = nicered]
 \tikzstyle{tbedge}=[ultra thick, draw=lightblue]
 \tikzstyle{redge}=[thick, draw = nicered]
 \tikzstyle{bedge}=[thick, draw = lightblue] 
 \tikzstyle{gedge}=[thick, draw = nicegreen] 
 \tikzstyle{brace} = [decorate, ultra thick, decoration = {calligraphic brace}]

\begin{document}

\maketitle

\begin{abstract}
    The {\sc Matching Cut} problem is to decide if the vertex set of a connected graph can be partitioned into two non-empty sets~$B$ and~$R$ such that the edges between~$B$ and~$R$ form a matching, that is, every vertex in~$B$ has at most one neighbour in~$R$, and vice versa.
    If for some integer $d\geq 1$, we allow every neighbour in~$B$ to have at most~$d$ neighbours in~$R$, and vice versa, we obtain the more general problem {\sc $d$-Cut}.
    It is known that {\sc $d$-Cut} is \NP-complete for every $d\geq 1$.
    However, for claw-free graphs, it is only known that {\sc $d$-Cut} is polynomial-time solvable for $d=1$ and \NP-complete for $d\geq 3$.
    We resolve the missing case $d=2$ by proving \NP-completeness.
    This follows from our more general study, in which we also bound the maximum degree.
    That is, we prove that for every  $d\geq 2$, {\sc $d$-Cut}, restricted to claw-free graphs of maximum degree~$p$, is constant-time solvable if $p\leq 2d+1$ and \NP-complete if $p\geq 2d+3$.
    Moreover, in the former case, we can find a $d$-cut in linear time.   
    We also show how our positive results for claw-free graphs can be generalized to $S_{1^t,\ell}$-free graphs where $S_{1^t,\ell}$ is the graph obtained from a star on $t+2$ vertices by subdividing one of its edges exactly $\ell$ times.
\end{abstract}

\section{Introduction}\label{s-intro}

In this paper, we determine new complexity results for finding $d$-cuts, which form a natural generalization of matching cuts.
The latter form a well-studied notion that combines two classic concepts in graph theory: matchings and edge cuts.
Namely, a {\it matching cut} in a connected graph~$G=(V,E)$ is a set of edges $M\subseteq E$ that is both a {\it matching} (i.e., no two edges in~$M$ share an end-vertex) and an {\it edge cut} (i.e., $V$ can be partitioned into two non-empty sets~$B$ and~$R$ such that~$M$ consists of all the edges between~$B$ and~$R$).
The {\sc Matching Cut} problem is to decide if a connected graph has a matching cut.

Already in 1984, Chv\'atal~\cite{Ch84} proved that {\sc Matching Cut} is \NP-complete even for $K_{1,4}$-free graphs of maximum degree~$4$.
Here, the graph $K_{1,\ell}$ denotes the {\it star} on $\ell+1$ vertices, that is, the graph with vertices $u,v_1,\ldots,v_\ell$ and edges~$uv_i$, for $i=1,\ldots,\ell$, and a graph is {\it $H$-free} if it does not contain~$H$ as an {\it induced} subgraph.
In contrast, Chv\'atal~\cite{Ch84} also showed that {\sc Matching Cut} is polynomial-time solvable for graphs of maximum degree at most~$3$.
In 2009, Bonsma~\cite{Bo09} proved the same for $K_{1,3}$-free graphs (of arbitrary degree), which are also known as {\it claw-free} graphs.
Since then a growing sequence of papers appeared in which the computational complexity of {\sc Matching Cut} for various special graph classes was determined, often in a systematic way (e.g.,~\cite{CHLLP21,FLPR23,LL23,LPR22,LPR23})  and also from a parameterized complexity perspective (e.g.~\cite{AKK22,GKKL22,GS21,KKL20,KL16}) and 
for many closely related problem variants, such as the problems of finding perfect matching cuts~\cite{BCD23,HT98,LT22}, maximum matching cuts~\cite{LPR24}, minimum matching cuts~\cite{LMO25}, matching multicuts~\cite{GJMS24} and disconnected perfect matchings~\cite{BP25}.
We refer to Section~\ref{s-con} for a summary of the complexity results for {\sc Matching Cut} restricted to $H$-free graphs as part of a summary for a more general notion, namely $d$-cuts, the topic of our paper.
For an integer~$d\geq 1$ and a connected graph $G=(V,E)$, a set $M\subseteq E$ is a {\it $d$-cut} of~$G$ if~$V$ can be partitioned into two non-empty sets~$B$ and~$R$ such that: 

\begin{itemize}
    \item the set~$M$ is the set of all edges between~$B$ and~$R$; and 
    \item every vertex in~$B$ has at most~$d$ neighbours in~$R$, and vice versa.
\end{itemize}

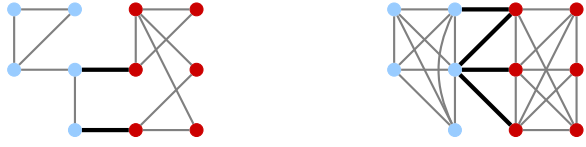
\begin{figure}[t]
    \centering
    \scalebox{1}{
\begin{tikzpicture}

\begin{scope}[scale = 0.8]
\node[bvertex] (v1) at (0,2){};
\node[bvertex] (v2) at (0,1){};
\node[bvertex] (v3) at (1,2){};
\node[bvertex] (v4) at (1,1){};
\node[bvertex] (v5) at (1,0){};
\node[rvertex] (u1) at (2,2){};
\node[rvertex] (u2) at (2,1){};
\node[rvertex] (u3) at (2,0){};
\node[rvertex] (u4) at (3,2){};
\node[rvertex] (u5) at (3,1){};
\node[rvertex] (u6) at (3,0){};

\draw[edge](v1) -- (v2);
\draw[edge](v1) -- (v3);
\draw[edge](v2) -- (v3);
\draw[edge](v2) -- (v4);
\draw[edge](v4) -- (v5);

\draw[tedge](v4) -- (u2);
\draw[tedge](v5) -- (u3);

\draw[edge](u1) -- (u2);
\draw[edge](u1) -- (u4);
\draw[edge](u1) -- (u5);
\draw[edge](u1) -- (u6);
\draw[edge](u2) -- (u4);
\draw[edge](u3) -- (u5);
\draw[edge](u3) -- (u6);

\end{scope}

\begin{scope}[shift = {(5,0)}, scale = 0.8]
\node[bvertex] (v1) at (0,2){};
\node[bvertex] (v2) at (0,1){};
\node[bvertex] (v3) at (1,2){};
\node[bvertex] (v4) at (1,1){};
\node[bvertex] (v5) at (1,0){};
\node[rvertex] (u1) at (2,2){};
\node[rvertex] (u2) at (2,1){};
\node[rvertex] (u3) at (2,0){};
\node[rvertex] (u4) at (3,2){};
\node[rvertex] (u5) at (3,1){};
\node[rvertex] (u6) at (3,0){};

\draw[edge] (v1) -- (v2);
\draw[edge] (v1) -- (v3);
\draw[edge] (v1) -- (v4);
\draw[edge] (v1) -- (v5);
\draw[edge] (v2) -- (v3);
\draw[edge] (v2) -- (v4);
\draw[edge] (v2) -- (v5);
\draw[edge] (v3) -- (v4);
\draw[edge] (v3) to [bend right = 25] (v5);
\draw[edge] (v4) -- (v5);

\draw[tedge] (v3) -- (u1);
\draw[tedge] (v4) -- (u1);
\draw[tedge] (v4) -- (u2);
\draw[tedge] (v4) -- (u3);

\draw[edge] (u1) -- (u2);
\draw[edge] (u1) -- (u4);
\draw[edge] (u1) -- (u6);
\draw[edge] (u2) -- (u3);
\draw[edge] (u2) -- (u4);
\draw[edge] (u2) -- (u5);
\draw[edge] (u2) -- (u6);
\draw[edge] (u3) -- (u4);
\draw[edge] (u3) -- (u5);
\draw[edge] (u3) -- (u6);
\draw[edge] (u4) -- (u5);
\draw[edge] (u5) -- (u6);

\end{scope}

\end{tikzpicture}
}
    \caption{Left: a graph with a matching cut ($1$-cut). Right: a graph with a $3$-cut but no $d$-cut for $d\leq 2$. Figure taken from~\cite{LMPS24}.}
    \label{f-examplesdcut}
\end{figure}

\noindent
See Figure~\ref{f-examplesdcut} for some examples.
We note that a $1$-cut is a matching cut and vice versa.
For a fixed integer~$d\geq 1$, the {\sc $d$-Cut} problem is to decide if a connected graph~$G$ has a $d$-cut.
Hence, {\sc Matching Cut} and {\sc $1$-Cut} are the same problems.

The {\sc $d$-Cut} problem was introduced by Gomes and Sau~\cite{GS21}.
Apart from several parameterized complexity results (see also \cite{AS21}), they
showed the following two results.
The first result is a structural result (shown for $d=1$ in~\cite{Ch84,Mo89}), which immediately implies the first statement of the second result (shown for $d=1$ in~\cite{Ch84}).

\begin{theorem}[Gomes and Sau~\cite{GS21}]\label{t-gs1}
    For $d\geq 1$, every graph $G=(V,E)$ with maximum degree $\Delta(G)\leq d+2$ and $|V|>7$ has a $d$-cut, which can be found in polynomial time.
\end{theorem}

\begin{theorem}[Gomes and Sau~\cite{GS21}]\label{t-gs2}
    For $d\geq 1$, {\sc $d$-Cut} is constant-time solvable for graphs with maximum degree $\Delta\leq d+2$ but \NP-complete for $(2d+2)$-regular graphs (that is, graphs in which every vertex has degree~$2d+2$).
\end{theorem}

\noindent
For $d\geq 1$, the results in Theorem~\ref{t-gs2} are currently the best results for {\sc $d$-Cut} of bounded maximum degree, so for $d\geq 2$ we have a complexity gap.
Gomes and Sau~\cite{GS21} argued that it is unlikely that the maximum degree bound of $2d+2$ in Theorem~\ref{t-gs2} can be lowered, as that would disprove a conjecture about the existence of so-called internal partitions of $r$-regular graphs for odd~$r$; see~\cite{BL16}.

In particular, Theorem~\ref{t-gs2} shows that {\sc $2$-Cut} is constant-time solvable for graphs of maximum degree $\Delta\leq 4$ and \NP-complete for graphs of maximum degree~$\Delta=6$.
Gomes and Sau~\cite{GS21} showed that all maximum degree-$5$ graphs on~$18$ vertices without a $2$-cut have some specific structure: these are $5$-regular or have exactly two vertices of degree~$4$, which are adjacent.
Interestingly, they mentioned in~\cite{GS21} that they were unable to find a maximum degree-$5$ graph on more than 18 vertices without a $2$-cut.

The {\sc $d$-Cut} problem has also been studied for $H$-free graphs.
As mentioned, we summarize these results in Section~\ref{s-con}.
For now, we focus on the case $H=K_{1,3}$ for the following reason.
It is known that for all $d\geq 3$, {\sc $d$-Cut} is \NP-complete for claw-free graphs. This was shown in~\cite{LMPS24}
(the gadget in the proof of Theorem~\ref{t-gs2} is not claw-free). 
The result for $d=3$ contrasts the case $d=1$ due to the aforementioned result of Bonsma~\cite{Bo09} that {\sc $1$-Cut} ({\sc Matching Cut}) is polynomial-time solvable on claw-free graphs and leaves open exactly the case $d=2$.
The authors of~\cite{LMPS24} could prove that {\sc $2$-Cut} is \NP-complete for $K_{1,4}$-free graphs and asked about the case $d=2$ in the following open problem:

\begin{open}[Lucke et al.~\cite{LMPS24}]\label{o-2}
    What is the complexity of {\sc $2$-Cut} for claw-free graphs?
\end{open}

\subsection*{Our Results}\label{s-ours}

Inspired by Theorems~\ref{t-gs1} and~\ref{t-gs2} of
Gomes and Sau~\cite{GS21} and motivated by Open Problem~\ref{o-2}, we consider {\it claw-free} graphs of bounded maximum degree.
In Sections~\ref{s-poly} and~\ref{s-np}, respectively, we show the following two results for $d\geq 2$, the second of which solves Open Problem~\ref{o-2}.

\begin{theorem}\label{t-existence}
    Let $d \geq 2$.
    Every \emph{claw-free} graph $G = (V,E)$ with maximum degree ${\Delta(G) \leq 2d+1}$ and ${|V|>4d^2(2d+1)}$ has a $d$-cut, which can be found in linear time.
    Moreover, there exist arbitrarily large claw-free $(2d+2)$-regular graphs with no $d$-cut.
\end{theorem}

\begin{theorem}\label{t-mainmain}
    For $d\geq 2$, {\sc $d$-Cut} is constant-time solvable for claw-free graphs with maximum degree $\Delta\leq 2d+1$ but \NP-complete for claw-free graphs with maximum degree $\Delta=2d+3$.
\end{theorem}

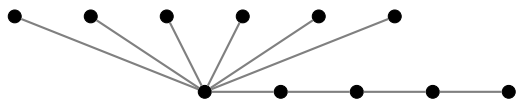
\begin{figure}
    \centering
    \begin{tikzpicture}

\node[vertex](u0) at (-1,0){};
\node[vertex](u1) at (0,0){};
\node[vertex](u2) at (1,0){};
\node[vertex](u3) at (2,0){};
\node[vertex](u4) at (3,0){};

\node[vertex](v1) at (-3.5,1){};
\node[vertex](v2) at (-2.5,1){};
\node[vertex](v3) at (-1.5,1){};
\node[vertex](v4) at (-0.5,1){};
\node[vertex](v5) at (0.5,1){};
\node[vertex](v6) at (1.5,1){};

\draw[edge](v1) -- (u0);
\draw[edge](v2) -- (u0);
\draw[edge](v3) -- (u0);
\draw[edge](v4) -- (u0);
\draw[edge](v5) -- (u0);
\draw[edge](v6) -- (u0);
\draw[edge](u0) -- (u1);
\draw[edge](u1) -- (u2);
\draw[edge](u2) -- (u3);
\draw[edge](u3) -- (u4);

\end{tikzpicture}
    \caption{The graph $S_{1^6,4}$.}
    \label{fig:broomstick}
\end{figure}

\noindent
The first statement of Theorem~\ref{t-existence} is the analogue of Theorem~\ref{t-gs1} for claw-free graphs 
(for which 
we can find the $d$-cut even in linear time).
We show how this statement follows from a corresponding stronger statement on $S_{1^t,\ell}$-free graphs, where $S_{1^t,\ell}$ is the graph obtained from a star on $t+2$ vertices by subdividing one of its edges exactly $\ell$ times, see Figure~\ref{fig:broomstick} for the case where $t=6$ and $\ell=4$.
The second statement of Theorem~\ref{t-existence} not only shows that the bound in the first statement is best possible, but also relates to Theorem~\ref{t-gs2}: it shows that for every $d\geq 2$, there exists an infinite number of $(2d+2)$-regular no-instances of {\sc $d$-Cut} that are also claw-free.

The first statement of Theorem~\ref{t-existence} immediately implies the constant-time part of Theorem~\ref{t-mainmain}, as $d$ is a constant.
We use the graph in the second statement of Theorem~\ref{t-existence} as part of our hardness gadget in the \NP-completeness part of Theorem~\ref{t-mainmain}.

Theorem~\ref{t-mainmain} shows in particular that {\sc $2$-Cut} is constant-time solvable for claw-free graphs of maximum degree $\Delta\leq 5$ but \NP-complete for claw-free graphs of maximum degree $\Delta=7$, which solves Open Problem~\ref{o-2}.
Moreover, the gadget used in the reduction for $d$-{\sc Cut} $(d\geq 3$) for claw-free graphs in~\cite{LMPS24} contains vertices of arbitrarily large degree.
Theorem~\ref{t-mainmain} strengthens this result by adding a maximum degree bound that is only one higher than the maximum degree bound in Theorem~\ref{t-gs2}.

\begin{figure}[t]
    \centering
    \begin{tikzpicture}

\begin{scope}[shift = {(0,0)}]
\node[vertex](v1) at (0,0){};
\node[vertex](v2) at (0.5,0.5){};
\node[vertex](v3) at (0.5,-0.5){};
\node[vertex](v4) at (1,0){};

\draw[edge](v1) -- (v2);
\draw[edge](v1) -- (v3);
\draw[edge](v2) -- (v3);
\draw[edge](v2) -- (v4);
\draw[edge](v3) -- (v4);
\end{scope}

\begin{scope}[shift = {(1,0)}]
\node[vertex](v1) at (0,0){};
\node[vertex](v2) at (0.5,0.5){};
\node[vertex](v3) at (0.5,-0.5){};
\node[vertex](v4) at (1,0){};

\draw[edge](v1) -- (v2);
\draw[edge](v1) -- (v3);
\draw[edge](v2) -- (v3);
\draw[edge](v2) -- (v4);
\draw[edge](v3) -- (v4);
\end{scope}

\begin{scope}[shift = {(2,0)}]
\node[vertex](v1) at (0,0){};
\node[vertex](v2) at (0.5,0.5){};
\node[vertex](v3) at (0.5,-0.5){};
\node[vertex](v4) at (1,0){};

\draw[edge](v1) -- (v2);
\draw[edge](v1) -- (v3);
\draw[edge](v2) -- (v3);
\draw[edge](v2) -- (v4);
\draw[edge](v3) -- (v4);
\end{scope}

\begin{scope}[shift = {(3,0)}]
\node[vertex](v1) at (0,0){};
\node[vertex](v2) at (0.5,0.5){};
\node[vertex](v3) at (0.5,-0.5){};
\node[vertex](v4) at (1,0){};

\draw[edge](v1) -- (v2);
\draw[edge](v1) -- (v3);
\draw[edge](v2) -- (v3);
\draw[edge](v2) -- (v4);
\draw[edge](v3) -- (v4);
\end{scope}

\end{tikzpicture}
    \caption{A chain of diamonds has no $1$-cut (the diamond is the graph $K_4-e$).}
    \label{fig:diamondchain}
\end{figure}
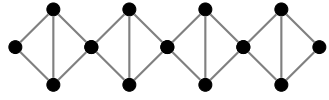

Note that~$d$ is at least~$2$ in both Theorems~\ref{t-existence} and~\ref{t-mainmain}, which can be explained as follows.
We first recall that for $d=1$, Bonsma~\cite{Bo09} proved that {\sc $d$-Cut} is polynomial time solvable even for claw-free 
graphs
with arbitrarily large maximum degree~$\Delta$.
Second, for $d=1$, Theorem~\ref{t-existence} would not hold: take, for example, a chain of graphs isomorphic to $K_p-e$ (complete graph minus an edge) which we obtain from $k$ copies of $K_p-e$, for $k \geq 1$, by identifying a pair of degree~$(p-2)$ vertices between consecutive copies.
That is, where $x^1_i$, $x^2_i$ is the pair of degree~$(p-2)$ vertices of the $i$th copy, we identify vertices $x^2_{i}$ and $x^1_{i+1}$, for all $1 \leq i \leq k-1$.
The resulting chain is an arbitrarily large claw-free graph of arbitrarily large maximum degree~$2p-4$ but with no $1$-cut, see Figure~\ref{fig:diamondchain} for the case where $p=4$.

In Section~\ref{s-con}, we give some directions for future work and
show the consequence of Theorem~\ref{t-mainmain} for the state-of-art summary of the complexity of {\sc $d$-Cut} for $H$-free graphs.

\section{Preliminaries}\label{l-pre}

In this paper, we only consider finite, undirected graphs without multiple edges and self-loops.
For every integer $n\geq1$, let $[n]:=\{1,\ldots,n\}$.

Let~$G = (V,E)$ be a graph and let~$v$ be a vertex of~$G$.
We denote by $N_G(v) = \{u \in V\; |\; uv \in E\}$ the \emph{(open) neighbourhood} of~$v$ and by $N_G[v] = N_G(v) \cup \{v\}$ the \emph{closed neighbourhood} of~$v$.
The \emph{degree} of~$v$ is the size of $N_G(v)$, and we denote by~$\Delta(G)$ the \emph{maximum degree} of~$G$.
For a set $S\subseteq V$, we denote by~$\delta(S)$ the set of edges of~$G$ with exactly one end-vertex in~$S$.

A {\it complete} graph is a graph whose vertex set is a {\it clique}, which is a set of pairwise adjacent vertices.
For $r\geq 1$, we let $K_r$ denote the complete graph on~$r$ vertices.
A vertex set $I \subseteq V$ is an \emph{independent set} of~$G$ if no two vertices of~$I$ are adjacent in~$G$.
In addition, $G$ is \emph{bipartite} if we can partition the vertex set into (possibly empty) sets~$A$ and~$B$ such that $A\cup B = V$ and~$A$ and~$B$ are each an independent set.

The {\it chromatic number} of a graph $G=(V,E)$ is the smallest integer~$k$ such that~$G$ has a {\it $k$-colouring}, which is a mapping $c:V\to \{1,\ldots,k\}$ with $c(u)\neq c(v)$ for every $uv\in E$. The set of all vertices that are mapped to the same colour $i\in \{1,\ldots,k\}$ is called a {\it colour class} of~$c$.
For ${k\geq 0}$, a graph~$G$ is \emph{$k$-degenerate} if every subgraph of~$G$ contains a vertex of degree at most~$k$. 
It is well known and readily seen that every $k$-degenerate graph is $(k+1)$-colourable.
The \emph{degeneracy} of~$G$ is the smallest integer~$k$ such that~$G$ is $k$-degenerate.
The following observation holds (we give a proof for completeness).

\begin{observation}\label{obs:mindeg}
    Let $t\geq2$ and let $G=(V,E)$ be a graph.
    If~$G$ has no independent set of size~$t$, then $G$ has a subgraph of minimum degree at least
    $\frac{|V|}{t-1}-1.$
\end{observation}

\begin{proof}
    Assume~$G$ has no independent set of size~$t$.
    Let~$k$ be the degeneracy of~$G$, so~$G$ has a subgraph of minimum degree at least~$k$.
    Consider a $(k+1)$-colouring~$c$ of~$G$. 
    As every colour class of~$c$ is an independent set, we have $k+1\geq \frac{|V|}{t-1}$ and thus $k\geq \frac{|V|}{t-1}-1$.
\end{proof}

\noindent
Let~$d \geq 1$ be an integer.
A \emph{red-blue colouring} of~$G$ gives each vertex of~$G$ either colour red or blue.
A \emph{red-blue $d$-colouring} of~$G$ is a red-blue colouring where every vertex has at most~$d$ neighbours of the other colour and both colours, red and blue, are used at least once.
See Figure~\ref{f-examplesdcut} for an example of a red-blue $1$-colouring (left figure) and a red-blue $3$-colouring (right figure).
In our proofs, we use the following well-known observation; see e.g.~\cite{LMPS24}.

\begin{observation}\label{obs-cut-colouring}
    For every integer $d \geq 1$, $G$ has a red-blue $d$-colouring with sets~$B$ and~$R$ if and only if the set of edges between~$B$ and~$R$ form a $d$-cut of~$G$.
\end{observation}

\noindent
Let~$G$ be a graph with a red-blue $d$-colouring for some $d\geq 1$.
We say that a set $S \subseteq V$ is \emph{monochromatic} if all vertices in~$S$ have the same colour.

We will also use the following observations, which can be readily seen.

\begin{observation}\label{l-safe}
    Let $d\geq 1$.
    Let $G=(V,E)$ be a graph and let $S\subseteq V$.
    If~$c_S$ is a red-blue colouring of $G[S]$ in which~$S$ is monochromatic, then every red-blue $d$-colouring of~$G$ that extends~$c_S$ gives each vertex with at least $d+1$ neighbours in~$S$ the same colour as the vertices of~$S$.
\end{observation}

\begin{observation}\label{l-clique}
    Let $d\geq 1$.
    Let~$G$ be a graph with a clique~$K$ on at least $2d+1$ vertices.
    Then~$K$ is monochromatic in every red-blue $d$-colouring of~$G$.
\end{observation}

\section{The Proof of Theorem~\ref{t-existence}}\label{s-poly}

In this section, we prove Theorem~\ref{t-existence}.
We start with the following lemma.
In its proof, we extend the argument
Chv\'atal~\cite{Ch84} used to observe that {\sc $1$-Cut} is 
constant-time solvable for graphs of maximum degree at most~$3$.

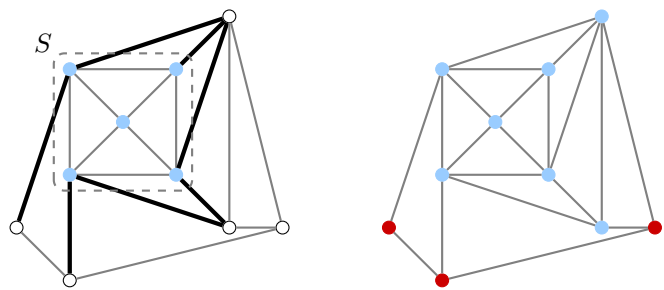
\begin{figure}
    \centering
    \begin{tikzpicture}
\begin{scope}[scale = 0.7]

    \node[bvertex](b1) at (1,2){};
    \node[bvertex](b2) at (1,4){};
    \node[bvertex](b3) at (3,4){};
    \node[bvertex](b4) at (3,2){};
    \node[bvertex](b5) at (2,3){};

    \node[evertex](b6) at (4,5){};
    \node[evertex](b7) at (4,1){};

    \node[evertex](r1) at (0,1){};
    \node[evertex](r2) at (1,0){};
    \node[evertex](r3) at (5,1){};

    \draw[edge](b1) -- (b2);
    \draw[edge](b3) -- (b2);
    \draw[edge](b3) -- (b4);
    \draw[edge](b1) -- (b4);

    \draw[edge](b1) -- (b5);
    \draw[edge](b2) -- (b5);
    \draw[edge](b3) -- (b5);
    \draw[edge](b4) -- (b5);

    \draw[tedge](b2) -- (b6);
    \draw[tedge](b3) -- (b6);
    \draw[tedge](b4) -- (b6);

    \draw[tedge](b1) -- (b7);
    \draw[tedge](b4) -- (b7);
    \draw[edge](b6) -- (b7);

    \draw[tedge](r1) -- (b2);
    \draw[tedge](r2) -- (b1);
    \draw[edge](r3) -- (b6);
    \draw[edge](r3) -- (b7);

    \draw[edge](r1) -- (r2);
    \draw[edge](r2) -- (r3);

    \node[](s) at ([shift = {(-0.5,0.5)}]b2) {$S$};
    \draw[edge, dashed, rounded corners=4pt] ([shift = {(-0.3,-0.3)}]b1) rectangle ([shift = {(0.3,0.3)}]b3);
\end{scope}

\begin{scope}[scale = 0.7, shift = {(7,0)}]

    \node[bvertex](b1) at (1,2){};
    \node[bvertex](b2) at (1,4){};
    \node[bvertex](b3) at (3,4){};
    \node[bvertex](b4) at (3,2){};
    \node[bvertex](b5) at (2,3){};

    \node[bvertex](b6) at (4,5){};
    \node[bvertex](b7) at (4,1){};

    \node[rvertex](r1) at (0,1){};
    \node[rvertex](r2) at (1,0){};
    \node[rvertex](r3) at (5,1){};

    \draw[edge](b1) -- (b2);
    \draw[edge](b3) -- (b2);
    \draw[edge](b3) -- (b4);
    \draw[edge](b1) -- (b4);

    \draw[edge](b1) -- (b5);
    \draw[edge](b2) -- (b5);
    \draw[edge](b3) -- (b5);
    \draw[edge](b4) -- (b5);

    \draw[edge](b2) -- (b6);
    \draw[edge](b3) -- (b6);
    \draw[edge](b4) -- (b6);

    \draw[edge](b1) -- (b7);
    \draw[edge](b4) -- (b7);
    \draw[edge](b6) -- (b7);

    \draw[edge](r1) -- (b2);
    \draw[edge](r2) -- (b1);
    \draw[edge](r3) -- (b6);
    \draw[edge](r3) -- (b7);

    \draw[edge](r1) -- (r2);
    \draw[edge](r2) -- (r3);
\end{scope}
\end{tikzpicture}   
    \caption{An example of the procedure in the proof of Lemma~\ref{lem:cubic} if ${d=2}$.
    Left: every vertex in~$S$ has at most~$2$ neighbours outside of~$S$. The edges in $\delta(S)$ are highlighted in black. Right: The red-blue $d$-colouring we obtain from the procedure.}
    \label{fig:deltaS}
\end{figure}

\begin{lemma}\label{lem:cubic}
    Let~$d \geq 1$.
    Let~$G = (V,E)$ be a graph with $\Delta(G) \leq 2d+1$, and let~$S$ be a non-empty set of vertices with ${|S|+|\delta(S)|<|V|}$.
    If each vertex in~$S$ is incident to at most~$d$ edges in~$\delta(S)$, then~$G$ has a $d$-cut, which can be found in linear time.
\end{lemma}
\begin{proof}
    Assume that each vertex in~$S$ is incident to at most~$d$ edges in~$\delta(S)$.
    We will now construct a red-blue $d$-colouring of~$G$, which,
    by Observation~\ref{obs-cut-colouring}, gives that $G$ has a $d$-cut.
    We first colour every vertex in~$S$ blue and recursively colour an uncoloured vertex blue if it has at least ${d+1}$ blue neighbours.
    Note that this takes $\mathcal{O}(|E|)$ time which is also linear in~$|V|$ as $|E|\leq(2d+1)|V|/2$.
    We colour all remaining vertices of~$G$ red.
    See also Figure~\ref{fig:deltaS}.
    
    Let~$R$ and~$B$ be the sets of red and blue vertices, respectively.
    Note that we found~$R$ and~$B$ in linear time and that $S \subseteq B$.
    By construction, every vertex in~$R$ has at most~$d$ neighbours in~$B$.
    As we assume that each vertex in~$S$ is incident to at most~$d$ edges in~$\delta(S)$, each vertex in~$S$ has at most~$d$ neighbours in~$R$.
    Since $\Delta(G) \leq 2d+1$ and each vertex in~$B\setminus S$ has at least $d+1$ neighbours in~$B$, each vertex in~$B\setminus S$ has at most~$d$ neighbours in~$R$.
    Hence, every vertex in~$B$ has at most~$d$ neighbours in~$R$.
    It remains to show that $R\neq \emptyset$, or equivalently, that $B\neq V$, which we will do below.

    Let~$B'$ be an arbitrary subset of~$B$ with $S\subseteq B'$.
    Recall that by construction, each vertex $v\in B'\setminus S$ has at least $d+1$ neighbours in~$B'$ and at most~$d$ neighbours not in~$B'$.
    So when we colour~$v$ blue, that is, add~$v$ to~$B$, we lose at least $d+1$ edges with exactly one blue end-vertex and gain at most~$d$ new edges with exactly one blue end-vertex.
    In other words, $|\delta(B'\setminus\{v\})|\geq|\delta(B')|+1$.
    This implies $|\delta(S)| \geq |\delta(B)| + (|B|- |S|)$, and thus $|\delta(S)| + |S| \geq |\delta(B)| + |B|$.
    We use this and our assumption that $|V| > |\delta(S)| + |S|$ to deduce that
    \[
        |V| > |\delta(S)| + |S| \geq |\delta(B)| + |B| \geq |B|,
    \]
    and thus~$B\neq V$.
    Hence, we obtained a red-blue $d$-colouring of~$G$.
\end{proof}

\noindent
To prove the first part of Theorem~\ref{t-existence}, we prove a stronger statement in Theorem~\ref{t-star-free} and afterwards we show how Theorem~\ref{t-star-free} implies the first part of Theorem~\ref{t-existence}.
For every positive integer~$t$ and~$\ell$, we recall that $S_{1^t,\ell}$ is the graph obtained from~$K_{1,t+1}$ by replacing one edge with a path of length~$\ell$, see Figure~\ref{fig:broomstick}.
Note that $S_{1^t,\ell}$ contains exactly ${t+\ell}$ edges.
In addition, $S_{1^t,1}$ is the graph $K_{1,t+1}$.

\begin{theorem}\label{t-star-free}
    For $d \geq 2$, $t\geq2$ and $\ell\geq 1$, every $S_{1^t,\ell}$-free graph $G = (V,E)$ with either maximum degree $\Delta=2$, or 
    $3\leq \Delta \leq \frac{t}{t-1}d+\frac{1}{t-1}$ and $|V|>(d+1)\left(\frac{\Delta(\Delta-1)^{\ell+1}-2}{\Delta-2}\right)$ has a $d$-cut, which can be found in linear time.
\end{theorem}

\begin{proof}   
    Let $d \geq 2$, $t\geq2$ and $\ell\geq 1$.
    Let~$G$ be an $S_{1^t,\ell}$-free graph.
If ${\Delta=2}$, then $G$ has a $d$-cut: take all edges between an any vertex of $G$ and its neighbours.
From now on, assume $3\leq \Delta \leq \frac{t}{t-1}d+\frac{1}{t-1}$ and and $|V|>(d+1)\left(\frac{\Delta(\Delta-1)^{\ell+1}-2}{\Delta-2}\right)$.

    We choose an arbitrary vertex~$v$ of~$G$ and let $S_0:=\{v\}$. 
    For each $i\in[\ell+1]$, we denote by~$S_i$ the set of vertices which are at distance exactly~$i$ from~$v$ in~$G$.
    Note that $|S_1|\leq\Delta$ and for each $i\in[\ell]$, $|S_{i+1}|\leq(\Delta-1)|S_i|$.
    Thus,
    \[
        \sum_{i=0}^{\ell+1}|S_i|\leq\frac{\Delta((\Delta-1)^{\ell+1}-1)}{\Delta-2}+1=\frac{\Delta(\Delta-1)^{\ell+1}-2}{\Delta-2}.
    \]
    
    Let~$U$ be the set of vertices in~$S_\ell$ which are incident to at least $d+1$ edges in $\delta(S_{\ell-1}\cup S_\ell)$.
    For each $u\in U$, we denote by~$N_u$ the set of neighbours of~$u$ not in $S_{\ell-1}\cup S_\ell$ and by~$G_u$ the subgraph of~$G$ induced 
    on~$N_u$.

    We first show that~$N_u$, for $u\in U$, has no independent set of size~$t$.
    Suppose not.
    By the definition of~$S_\ell$, there is an induced path of length~$\ell$ between~$v$ and~$u$.
    Together with this path and the independent set in~$N_u$, we can find $S_{1^t,\ell}$ as an induced subgraph, a contradiction.
    Thus, $N_u$ has no independent set of size~$t$.

    For each $u\in U$, since~$G_u$ has no independent set of size~$t$, by Observation~\ref{obs:mindeg}, $G_u$ has a subgraph~$H_u$ whose minimum degree is at least
    \[
        \frac{|V(G_u)|}{t-1}-1\geq\frac{d+1}{t-1}-1.
    \]
    We remark that~$H_u$ can be found in constant time as $|N_u|\leq\Delta-1$, which is a constant.
    Now, we let
    \[
        S:=\left(\bigcup_{i=0}^\ell S_i\right)\cup\left(\bigcup_{u\in U}V(H_u)\right).
    \]
    Note that~$S$ can be found in constant time as $|U|\leq|S_\ell|\leq\Delta(\Delta-1)^{\ell-1}$, which is a constant.

    We are going to apply Lemma~\ref{lem:cubic} to~$S$.
    Since $t\geq2$, we have $t\leq2(t-1)$, and therefore
    \[
        \Delta\leq\frac{t}{t-1}d+\frac{1}{t-1}\leq2d+1.
    \]
    To apply Lemma~\ref{lem:cubic} to~$S$, we first show that every vertex in~$S$ is incident to at most~$d$ edges in~$\delta(S)$.
    By construction, this holds for every vertex in $S\setminus\bigcup_{u\in U}(N_u\cup\{u\})$.
    Let~$u$ be a vertex in~$U$ and let~$w$ be a vertex in $N_u\cup\{u\}$.
    Note that~$u$ and~$w$ might be the same.
    Since the minimum degree of~$H_u$ is at least
    \[
        \frac{d+1}{t-1}-1
    \]
    and~$u$ is adjacent to every vertex in~$N_u$, there are at least
    \[
        \frac{d+1}{t-1}
    \]
    neighbours of~$w$ in~$S$.
    Then~$w$ is incident to at most
    \[
        \Delta-\frac{d+1}{t-1}
        \leq\frac{t}{t-1}d+\frac{1}{t-1}-\frac{d+1}{t-1}=d
    \]
    edges in~$\delta(S)$.
    That is, $w$ is incident to at most~$d$ edges in~$\delta(S)$.
    Hence, every vertex in~$S$ is incident to at most~$d$ edges in~$\delta(S)$.

    We now show that $|V|>|S|+|\delta(S)|$.
    Since~$S\subseteq\bigcup_{i=0}^{\ell+1}S_i$, we have
    \[
        |S|\leq\sum_{i=0}^{\ell+1}|S_i|\leq\frac{\Delta(\Delta-1)^{\ell+1}-2}{\Delta-2}.
    \]
    Since every vertex in~$S$ is incident to at most~$d$ edges in~$\delta(S)$, we have $|\delta(S)|\leq d|S|$.
    Thus,
    \[
        |V|>(d+1)\left(\frac{\Delta(\Delta-1)^{\ell+1}-2}{\Delta-2}\right)\geq(d+1)|S|=|S|+d|S|\geq|S|+|\delta(S)|.
    \]
    Hence, by Lemma~\ref{lem:cubic}, we can find a $d$-cut of~$G$ in linear time.
\end{proof}

\noindent
We now show that Theorem~\ref{t-star-free} implies Theorem~\ref{t-existence}, which we recall below.

\medskip
\noindent
{\bf Theorem~\ref{t-existence} (first part, restated).}
{\it For $d \geq 2$, every claw-free graph $G = (V,E)$ with maximum degree ${\Delta(G) \leq 2d+1}$ and ${|V|>4d^2(2d+1)}$ has a $d$-cut, which can be found in linear time.}

\begin{proof}
Let $G$ be a claw-free graph.
If $\Delta(G)=2$, then we can apply 
Theorem~\ref{t-star-free} directly.
Assume $\Delta(G)\geq 3$.
Theorem~\ref{t-star-free} with $t=2$ and $\ell=1$ implies the first part of Theorem~\ref{t-existence}.
In order to see this, let $G$ be a claw-free graph with maximum degree ${\Delta(G) \leq 2d+1}$ and ${|V|>4d^2(2d+1)}$.
As $2d+1=\frac{t}{t-1}d+\frac{1}{t-1}$ if $t=2$, we obtain $\Delta(G)\leq\frac{t}{t-1}d+\frac{1}{t-1}$.
We also observe:
\[
    |V|>4d^2(2d+1)=(2d-1)\left(\frac{(2d+1)(2d)^2}{2d-1}\right)>(d+1)\left(\frac{\Delta(\Delta-1)^2-2}{\Delta-2}\right)
\]
where the last inequality holds from the fact that $2d-1\geq d+1$ and the function ${x(x-1)^2/(x-2)}$ is increasing for $x\geq 3$. Hence,
the conditions in Theorem~\ref{t-star-free} are satisfied, and we may conclude that $G$ has a $d$-cut, which can be found in linear time.
\end{proof}

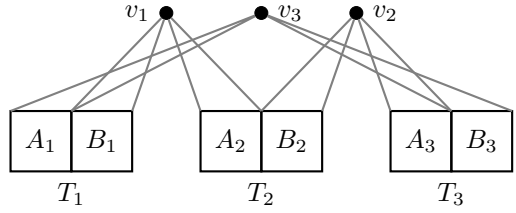
\begin{figure}
    \centering



\begin{tikzpicture}
\tikzstyle{clique}=[thick,rectangle,inner sep=0.cm, minimum size=8mm, fill=none, draw=black]

    \node[clique](c1) at (0,0){$A_1$};
    \node[clique](d1) at (0.8,0){$B_1$};

    \node[clique](c2) at (2.5,0){$A_2$};
    \node[clique](d2) at (3.3,0){$B_2$};

    \node[clique](c3) at (5,0){$A_3$};
    \node[clique](d3) at (5.8,0){$B_3$};

    \node[vertex,label=left: $v_1$](v1) at (1.65,1.7){};
    \draw[edge] (0.4,0.4) -- (v1) -- (1.2,0.4);
    \draw[edge] (2.1,0.4) -- (v1) -- (2.9,0.4);

    \node[vertex,label=right: $v_2$](v2) at (4.15,1.7){};
    \draw[edge] (2.9,0.4) -- (v2) -- (3.7,0.4);
    \draw[edge] (4.6,0.4) -- (v2) -- (5.4,0.4);

    \node[vertex,label=right: $v_3$](v3) at (2.9,1.7){};
    \draw[edge] (-0.4,0.4) -- (v3) -- (0.4,0.4);
    \draw[edge] (5.4,0.4) -- (v3) -- (6.2,0.4);

    \node[] (t1) at (0.4, -0.7){$T_1$};
    \node[] (t2) at (2.9, -0.7){$T_2$};
    \node[] (t3) at (5.4, -0.7){$T_3$};

\end{tikzpicture}
    \caption{The graph constructed in Theorem~\ref{thm:examples} for ${k=3}$. Note that for each $i\in[3]$, ${A_i\cup B_i}$ is a clique and~$v_i$ is complete to $B_i\cup A_{i+1}$, where $A_4:=A_1$,}
    \label{fig:graph-wo-2cut}
\end{figure}

\noindent
We now prove the second part of Theorem~\ref{t-existence}, which shows that the bound given in the first part of Theorem~\ref{t-existence} is best possible.
We slightly reformulate the statement as follows.

\begin{theorem}\label{thm:examples}
    For every integer $d\geq2$, $k\geq2$, and $r\geq 2d+2$, there exists a claw-free $r$-regular graph on $(r+1)k$ vertices that has no $d$-cut.
\end{theorem}
\begin{proof}
    Let $T_1,\ldots,T_k$ be pairwise disjoint cliques of size~$r$.
    For each $i\in[k]$, let $\{A_i,B_i\}$ be a partition of~$T_i$ such that~${|A_i|=d+1}$.
    Note that
    \[
        |B_i|=r-(d+1)\geq d+1.
    \]
    Let~$H$ be the graph obtained from $T_1\cup\cdots\cup T_k$ by adding~$k$ new vertices $v_1,\ldots,v_k$ such that for each $i\in[k]$, $N_H(v_i)$ is equal to $B_i\cup A_{i+1}$ where ${A_{k+1}:=A_1}$, see Figure~\ref{fig:graph-wo-2cut}.
    Note that~$H$ has ${rk+k=(r+1)k}$ vertices.

    We first show that~$H$ is claw-free and $r$-regular.
    Let~$i$ be an arbitrary integer in~$[k]$.
    Since the neighbourhood of~$v_i$ is the disjoint union of two cliques~$B_i$ and~$A_{i+1}$, $H$ has no claw centred at~$v_i$.
    Note further that $v_i$ has degree~$r$.
    Consider now a vertex $v$ in~$T_i$.
    The neighbourhood of $v$ consists of all other vertices in $T_i$ and some vertex $v_j$, where $j \in \{i-1,i,k\}$. 
    Hence, $v$ has degree ${(r-1)+1=r}$ and since $T_i$ is a clique, $H$ has no claw centred at~$v$.
    We conclude that~$H$ is claw-free and $r$-regular.

    It remains to show that~$H$ has no $d$-cut.
    Note that in any red-blue $d$-colouring of $H$, $T_i$, for $i\in[k]$, is monochromatic by Observation~\ref{l-clique}, since it is a clique of size at least $2d+2$.
    We claim that  
    \[
        T:=\bigcup_{i=1}^k T_i
    \]
    is monochromatic in every red-blue $d$-colouring of $H$.
    Suppose for a contradiction that $H$ has a red-blue $d$-colouring where $T$ is not monochromatic.
    Then there exists $j\in[k]$ such that~$T_j$ and~$T_{j+1}$ have different colours, where ${T_{k+1}:=T_1}$.
    Without loss of generality, assume that~$T_j$ is coloured red.
    Consequently, $T_{j+1}$ is coloured blue.
    Since~$v_j$ has at least $d+1$ red neighbours in~$B_j$, it is coloured red.
    However, $v_j$ has $d+1$ blue neighbours in~$A_{j+1}$, contradicting the fact that we considered a red-blue $d$-colouring.
    Hence, if~$H$ has a red-blue $d$-colouring, then~$T$ is monochromatic.

    Thus, if~$H$ has a red-blue $d$-colouring, then $T$ is monochromatic and hence, $v_i$, for $i\in[k]$, is coloured the same as~$T$ since~$v_i$ has at least $2d+2$ neighbours in~$T$.
    Therefore, there is no red-blue $d$-colouring of~$H$.
    By Observation~\ref{obs-cut-colouring}, $H$ has no $d$-cut.
\end{proof}

\section{The Proof of Theorem~\ref{t-mainmain}}\label{s-np}

In this section we prove Theorem~\ref{t-mainmain}.
We first recall that
the first statement of Theorem~\ref{t-existence} immediately implies the constant-time part of Theorem~\ref{t-mainmain}, as $d$ is a constant.
Hence, it remains to show the \NP-completeness part of Theorem~\ref{t-mainmain}, that is, that for every integer $d \geq 2$, {\sc $d$-Cut} is \NP-complete for claw-free graphs of maximum degree ${2d+3}$.
To do this, we will reduce from {\sc NAE $3$-Sat $0$-$1$}, which is a variant of {\sc Not-All-Equal $3$-Satisfiability (NAE $3$-Sat)}.
An instance of {\sc NAE $3$-Sat} consists of a set of clauses $C_1,\ldots,C_m$, each containing three variables.
The problem is to decide whether there is an assignment of the variables such that every clause contains both a true and a false literal.
Such an assignment is called an \emph{NAE-satisfying} assignment.

An instance of the problem we will reduce from {\sc NAE $3$-Sat $0$-$1$} is an instance of {\sc NAE $3$-Sat} that is satisfied both by the assignment where every variable is true and the assignment where every variable is false.
In other words, given an instance of {\sc NAE $3$-Sat $0$-$1$}, every clause contains both of a positive literal and a negative literal.
Additionally, since the clause $(\bar{x} \lor \bar{y} \lor z)$ has the same set of NAE-satisfying assignments as the clause $(x \lor y \lor \bar{z})$, we may assume that each clause contains exactly one negative literal.

The problem {\sc NAE $3$-Sat $0$-$1$} is to decide whether an instance has a third NAE-satisfying assignment.
That is, we must decide whether there exists an assignment of variables such that at least one variable is true, at least one variable is false, and each clause contains both a true and a false literal.

\begin{proposition}[Eagling-Vose et al.~\cite{EMPS25}]\label{prop:EMPS}
{\sc NAE $3$-Sat $0$-$1$} is \NP-complete. 
\end{proposition}

\noindent
We remark that it is possible to combine the reduction from {\sc Matching Cut} to {\sc NAE $3$-Sat $0$-$1$} given in \cite{EMPS25} with the reduction given here from the latter problem to {\sc $d$-Cut} to obtain a reduction directly from {\sc Matching Cut} to {\sc $d$-Cut}.
However, we separate the two reductions for ease of explanation.

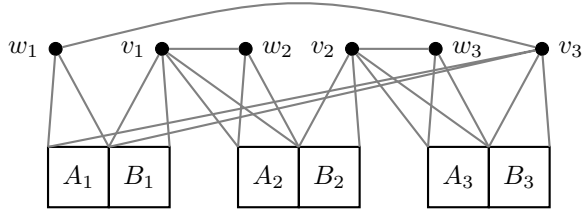
\begin{figure}
    \centering
    \begin{tikzpicture}
\tikzstyle{clique}=[thick,rectangle,inner sep=0.cm, minimum size=8mm, fill=none, draw=black]

    \node[clique](c1) at (0,0){$A_1$};
    \node[clique](d1) at (0.8,0){$B_1$};

    \node[clique](c2) at (2.5,0){$A_2$};
    \node[clique](d2) at (3.3,0){$B_2$};

    \node[clique](c3) at (5,0){$A_3$};
    \node[clique](d3) at (5.8,0){$B_3$};

    \node[vertex,label=left: $w_1$](w1) at (-0.3,1.7){};
    \draw[edge] (-0.4,0.4) -- (w1) -- (0.4,0.4);

    \node[vertex,label=left: $v_1$](v1) at (1.1,1.7){};
    \draw[edge] (0.4,0.4) -- (v1) -- (1.2,0.4);
    \draw[edge] (2.1,0.4) -- (v1) -- (2.9,0.4);

    \node[vertex,label=right: $w_2$](w2) at (2.2,1.7){};
    \draw[edge] (2.1,0.4) -- (w2) -- (2.9,0.4);
    \draw[edge] (v1) -- (w2);

    \node[vertex,label=left: $v_2$](v2) at (3.6,1.7){};
    \draw[edge] (2.9,0.4) -- (v2) -- (3.7,0.4);
    \draw[edge] (4.6,0.4) -- (v2) -- (5.4,0.4);

    \node[vertex,label=right: $w_3$](w3) at (4.7,1.7){};
    \draw[edge] (4.6,0.4) -- (w3) -- (5.4,0.4);
    \draw[edge] (v2) -- (w3);

    \node[vertex,label=right: $v_3$](v3) at (6.1,1.7){};
    \draw[edge] (-0.4,0.4) -- (v3) -- (0.4,0.4);
    \draw[edge] (5.4,0.4) -- (v3) -- (6.2,0.4);
    \draw[edge] (v3) .. controls (2.9,2.5) .. (w1);
\end{tikzpicture}
    \caption{The graph illustrating $H_{d,k,r}$ for ${k=3}$. Note that for each $i\in[3]$, ${A_i\cup B_i}$ is a clique, $v_i$ is complete to $B_i\cup A_{i+1}$, where $A_4:=A_1$, and~$w_i$ is complete to~$A_i$. If we remove~$w_1$, $w_2$, and~$w_3$, we obtain the graph in Figure~\ref{fig:graph-wo-2cut}.}
    \label{fig:graph-Hdkr}
\end{figure}

For the reduction in the proof of our \NP-hardness result, we use a lemma that is based on the class of graphs from Theorem~\ref{thm:examples}.
For every integer $d\geq2$, $k\geq2$, and $r\geq 2d+2$, let $H_{d,k,r}$ be the graph obtained from the graph in Theorem~\ref{thm:examples}, with respect to the same integers, by adding~$k$ new vertices $w_1,\ldots,w_k$ such that for each $i\in[k]$, the neighbourhood of~$w_i$ in~$H_{d,k,r}$ is equal to $A_i\cup\{v_{i-1}\}$ where $v_0:=v_k$.
See also Figure~\ref{fig:graph-Hdkr}.

\begin{lemma}\label{cor:examples}
    For every integer $d\geq2$, $k\geq2$, and $r\geq 2d+2$, the graph~$H_{d,k,r}$ is a claw-free graph with maximum degree $r+1$ on $(r+2)k$ vertices which has no $d$-cut.
    In particular, each of $w_1,\ldots,w_k$ has degree $d+2$, and the other vertices have degree at least~$r$.
\end{lemma}

\begin{proof}
    Let~$H$ be the graph in Theorem~\ref{thm:examples}, with respect to the same integers, and let $H':=H_{d,k,r}$.
    Note that
    \[
        |V(H')|=|V(H)|+k=(r+1)k+k=(r+2)k.
    \]
    Since~$H$ is $r$-regular and $w_1,\ldots,w_k$ have pairwise disjoint neighbourhoods, each of $w_1,\ldots,w_k$ has degree $d+2$, and the other vertices have degree either~$r$ or~$r+1$.

    We show that~$H'$ is claw-free.
    Suppose that~$H'$ has a claw~$K$ centred at a vertex~$c$.
    Since~$H$ is claw-free, $K$ contains~$w_i$ for some $i\in[k]$.
    As $N_{H'}(w_i)$ is a clique, $c\neq w_i$.
    So $N_{H'}(c)$ is the union of two cliques, contradicting that~$K$ is a claw.
    Hence, $H'$ is claw-free.

    It remains to show that~$H'$ has no $d$-cut.
    By Observation~\ref{obs-cut-colouring}, $H'$ has a $d$-cut if and only if~$H'$ has a red-blue $d$-colouring.
    By Theorem~\ref{thm:examples}, we know that $V(H)$ is monochromatic in every red-blue $d$-colouring.
    Further, $w_i$, for $i\in[k]$, has $d+2$ neighbours in~$H$.
    Thus, by Observation~\ref{l-safe}, $w_i$ is coloured the same as~$V(H)$.
    This implies that $V(H')$ is monochromatic in every red-blue $d$-colouring and hence $H'$ has no red-blue $d$-colouring.
    It follows that~$H'$ has no $d$-cut.
\end{proof}

\noindent
We are now ready to prove the remaining part of Theorem~\ref{t-mainmain}, which we restate below.

\medskip
\noindent
{\bf Theorem~\ref{t-mainmain} (\NP-completeness part, restated).}
{\it For $d\geq 2$, {\sc $d$-Cut} is \NP-complete for claw-free graphs with maximum degree $\Delta=2d+3$.}

\begin{proof}
    We reduce from {\sc NAE $3$-Sat $0$-$1$}.
    Let~$\phi$ be an instance of {\sc NAE $3$-Sat $0$-$1$} with variables $x_1,\dots, x_n$ and clauses $C_1, \dots, C_m$.
    In the following, we construct an instance~$G$ of {\sc $d$-Cut} which is claw-free and has maximum degree $2d+3$. 
    For each $\ell\in[n]$, let~$k_\ell$ be the number of occurrences of~$x_\ell$ in~$\phi$ and let~$F_\ell$ be a copy of~$H_{d,k_\ell,2d+2}$.
    We say that a vertex of~$F_\ell$ is \emph{free} if it has degree $d+2$ in~$F_\ell$.
    Note that~$F_\ell$ has exactly~$k_\ell$ free vertices.

    For each clause $C_i=(\overline{x}_{i_1} \lor x_{i_2} \lor x_{i_3} )$, for $i \in [m]$, we add two disjoint cliques~$D_{i,1}$ of size~$d$ and~$D_{i,2}$ of size $d+1$.
    In addition, we add a vertex $c_i$.
    We let $D_i:=D_{i,1}\cup D_{i,2}$ and call $D_i\cup \{c_i\}$ a \emph{clause gadget}.
    We add all edges between~$D_{i,1}$ and~$D_{i,2}$ and between~$c_i$ and~$D_{i,1}$.
    For each $j\in[3]$, we choose some free vertex of~$F_{i_j}$ and call this vertex~$w_{i_j}$.
    We add edges between~$w_{i_1}$ and each vertex in~$D_{i,1}\cup\{c_i\}$, edges between~$w_{i_2}$ and each vertex in~$D_{i,2}$, and we also add an edge between~$c_i$ and~$w_{i_3}$, see Figure~\ref{fig:hard-claw}.
    We choose these free vertices in such a way that every free vertex has neighbours in exactly one clause gadget.
    This is possible as, for every $\ell\in[n]$, $F_\ell$ has~$k_\ell$ free vertices and~$x_\ell$ appears~$k_\ell$ times in~$\phi$.
    Let~$G$ be the resulting graph.

\begin{figure}[t]
    \centering
    \begin{tikzpicture}

\tikzstyle{clique}=[thick,rectangle,inner sep=0.cm, minimum height=10mm, minimum width = 15mm, fill=none, draw=black]
\tikzstyle{bigclique}=[thick,rectangle,inner sep=0.cm, minimum size=15mm, fill=none, draw=black]

    \node[clique](box1) at (2,1){};
    \node[vertex,label=below: $w_{i_1}$](wi1) at (box1.north){};

    \node[label = right: $F_{i_1}$](t) at (box1.east){};

    \node[clique](box2) at (5.0,1){};
    \node[vertex,label=below: $w_{i_2}$](wi2) at (box2.north){};
    
    \node[label = right: $F_{i_2}$](t) at (box2.east){};

    \node[clique](box3) at (0,1){};
    \node[vertex,label=below: $w_{i_3}$](wi3) at (box3.north){};
    
     \node[label = left: $F_{i_3}$](t) at (box3.west){};

\node[bigclique] (Di1) at (2.75,3.5){$D_{i,1}$};
\node[bigclique] (Di2) at (4.25,3.5){$D_{i,2}$};
\node[vertex, label = left: $c_i$] (ci) at (0.5,3){};

\draw[edge] (ci) -- (Di1.north west);
\draw[edge] (ci) -- (Di1.south west);

\draw[edge] (ci) -- (wi3);
\draw[edge] (wi1) -- (ci);
\draw[edge] (wi1) -- (Di1.south east);
\draw[edge] (wi1) -- (Di1.south west);

\draw[edge] (wi2) -- (Di2.south east);
\draw[edge] (wi2) -- (Di2.south west);

\end{tikzpicture}
    \caption{The clause gadget consisting of the vertex $c_{i}$ and the clique $D_{i,1}\cup D_{i,2}$ representing the clause $C_i=(\overline{x}_{i_1} \lor x_{i_2} \lor x_{i_3} )$, together with the corresponding variable vertices.}
    \label{fig:hard-claw}
\end{figure}
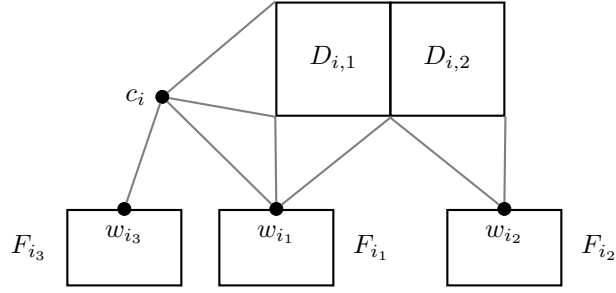

    \medskip
    We first show that~$G$ is a claw-free graph of maximum degree $2d+3$.
    Let~$\ell$ be an arbitrary integer in~$[n]$.
    Let $v \in V(F_\ell)$ be a vertex that is not free.
    Then, by Lemma~\ref{cor:examples}, $v$ has degree at most~$2d+3$ and further, since $v$ has neighbours only in~$F_\ell$, $G$ has no claw centred at~$v$.
    Let~$w$ be a free vertex of~$F_\ell$.
    Recall that~$w$ has a neighbour in a clause gadget~$D_i \cup \{c_i\}$ for some $i\in[m]$.
    Note that~$w$ has degree at most $(d+2)+(d+1)=2d+3$ since it has $d+2$ neighbours within~$F_\ell$ and at most $d+1$ neighbours in~$D_i\cup \{c_i\}$.
    In addition, the neighbours of~$w$ in~$D_i\cup \{c_i\}$ form a clique of size~$1$ or~$d+1$.
    Thus, the neighbourhood of~$w$ in~$G$ is the union of this clique in $D_i\cup \{c_i\}$ and some clique in~$F_\ell$, so~$G$ has no claw centred at~$w$.
    
    Finally, we consider the clause gadget~$D_i \cup \{c_i\}$.
    Recall that $c_i$ is a fixed vertex with neighbourhood $D_{i,1}\cup\{w_{i_1},w_{i_3}\}$.
    Thus, it has degree ${d+2}$.
    Since $w_{i_1}$ is complete to $D_{i,1}$, the neighbourhood of $c_i$ is the union of two cliques.
    Hence, the vertex~$c_i$ is not the centre of a claw of~$G$.
    Let~$u$ be a vertex in $D_{i,1}$.
    The closed neighbourhood of~$u$ is $D_i\cup\{c_i, w_{i_1}\}$, which is the union of the two cliques $D_{i,1}\cup \{c_i, w_{i_1}\}$ and $D_{i_2}$.
    Hence, $u$ is not the centre of a claw in~$G$.
    Since ${|D_i|=2d+1}$, $u$ has degree ${2d + 2}$.
    Now, let~$u'$ be a vertex in~$D_{i,2}$.
    The closed neighbourhood of~$u'$ is $D_i\cup\{w_{i_2}\}$, which is the union of the two cliques $D_i$ and~$\{w_{i_2}\}$.
    So~$u'$ is not the centre of a claw in~$G$.
    Again, as ${|D_{i}|=2d+1}$, the vertex~$u'$ has degree $2d+1$.
    It follows that~$G$ is a claw-free graph of maximum degree at most $2d+3$.

    \medskip
    We show that~$G$ has a $d$-cut if and only if~$\phi$ has an NAE-satisfying assignment with both a true and a false variable.
    By Observation~\ref{obs-cut-colouring}, this is equivalent to showing that~$G$ has a red-blue $d$-colouring if and only if~$\phi$ has an NAE-satisfying assignment with both a true and a false variable.
    
    We assume first that~$G$ has a red-blue $d$-colouring and show that $\phi$ has an NAE-satisfying assignment with both a true and a false variable.
    Recall that by Lemma~\ref{cor:examples}, $F_\ell$ is monochromatic for each $\ell\in[n]$.
    We set a variable~$x_\ell$ to true if $F_\ell$ is coloured blue and to false otherwise.

    We first show that the resulting assignment contains both a true and a false variable.
    Suppose for a contradiction that this is not the case.
    We may assume without loss of generality that all variables are true.
    This implies that for every $\ell\in[n]$, $F_\ell$ is coloured blue.
    Note that $D_i$, for $i\in[m]$, is a clique of size $2d+1$ and thus monochromatic by Observation~\ref{l-clique}.
    As $w_{i_2} \in F_{i_2}$, it is coloured blue.
    Further, $w_{i_2}$ has $d+1$ neighbours in~$D_{i,2} \subseteq D_i$.
    Hence, since $D_{i,2}$ is monochromatic, $D_{i,2}$ is coloured blue.
    It follows that $D_{i,1}$ is blue as well and, since $c_i$ has $d+1$ blue neighbours in $D_{i,1}\cup\{w_{i_1}\}$, it is coloured blue.
    Thus, the whole graph is coloured blue, contradicting the fact that we considered a red-blue $d$-colouring.
    Therefore, the assignment contains both a true and a false variable.

    \smallskip
    We now show that the assignment is NAE-satisfying for~$\phi$.
    Suppose for a contradiction that the assignment is not NAE-satisfying.
    Then there exists a clause $C_i=(\overline{x}_{i_1} \lor x_{i_2} \lor x_{i_3})$, for some $i \in [m]$, such that all literals have the same value.
    Note that this implies that both~$x_{i_2}$ and~$x_{i_3}$ take the opposite value of~$x_{i_1}$.
    Without loss of generality, assume that $x_{i_1}$ is true and $x_{i_2}$ and~$x_{i_3}$ are both false.
    That is, $w_{i_1}$ is coloured blue while~$w_{i_2}$ and~$w_{i_3}$ are coloured red.
    Recall that $D_i$ is monochromatic.
    Since $w_{i_2}$ has $d+1$ neighbours in~$D_{i,2}$, $D_i$ is coloured the same as $w_{i_2}$ and is thus coloured red.
    This implies that the blue vertex $w_{i_1}$ has~$d$ red neighbours in~$D_{i,1}$, so its neighbour~$c_i$ is coloured blue.
    However, as~$w_{i_3}$ is coloured red, $c_i$ has $d+1$ red neighbours in $D_{i,1}\cup\{w_{i_3}\}$, contradicting the fact that we considered a red-blue $d$-colouring.
    Hence, the assignment is NAE-satisfying for~$\phi$.

    \medskip
    For the other direction, we now assume that~$\phi$ has an NAE-satisfying assignment with both a true and a false variable.
    We construct a red-blue $d$-colouring of~$G$ as follows.
    For each $\ell\in[n]$, we colour~$F_\ell$ blue if it corresponds to a true variable and red otherwise.
    For each $i\in[m]$, we colour $D_i$ the same as~$w_{i_2}$.
    In addition, we colour~$c_i$ the same as~$w_{i_1}$.

    We show that this colouring is indeed a red-blue $d$-colouring of $G$.
    Note first that since the assignment contains both a true and a false variable we have both red and blue vertices.
    It remains to show that every vertex has at most~$d$ neighbours of the other colour.
    Let~$\ell$ be an arbitrary integer in~$[n]$.
    Note that every vertex of~$F_\ell$ which is not a free vertex has no neighbour of the other colour.
    Let~$w$ be a free vertex of~$F_\ell$.
    Recall that there is exactly one clause $C_i$, for $i\in[m]$, such that $w$ is adjacent to the corresponding clause gadget $D_i$.
    That is, $w\in\{w_{i_1},w_{i_2},w_{i_3}\}$.
    We distinguish between these cases:
    
    \begin{itemize}
        \item If $w=w_{i_3}$, then it has one neighbour outside of $F_\ell$ and thus at most one neighbour of the other colour.
        \item If $w=w_{i_2}$, then, by construction, it has no neighbour of the other colour.
        \item If $w=w_{i_1}$, then it has $d+1$ neighbours outside of~$F_\ell$ and, by construction, at least one of them, namely $c_i$, has the same colour as~$w$.
    \end{itemize}

\noindent    
    Thus, $w$ has at most~$d$ neighbours of the other colour.
    
    Now, let~$u$ be a vertex in~$D_i$.
    If $u\in D_{i,2}$, then its closed neighbourhood is $D_i\cup\{w_{i_2}\}$, which is monochromatic by construction.
    Thus, $u$ has no neighbour of the other colour.
    If $u\in D_{i,1}$, then, since $D_i$ is monochromatic, the only neighbours of $u$ which may be coloured with the other colour are $w_{i_1}$ and~$c_i$.
    Since $d\geq2$, it follows that $u$ has at most~$d$ neighbours of the other colour.
    Thus, we may assume that $u=c_i$.
    Recall that the neighbourhood of~$c_i$ is $D_{i,1}\cup\{w_{i_1},w_{i_3}\}$ and that~$c_i$ is coloured the same as~$w_{i_1}$.
    Assume without loss of generality that $c_i$ is coloured blue.
    Suppose for a contradiction that $c_i$ has at least $d+1$ red neighbours.
    That is, $D_{i,1}$ and $w_{i,3}$ are coloured red.
    Since $D_i$ is monochromatic by construction, this implies that $D_{i,2}$ and thus $w_{i,2}$ are coloured red as well.
    Hence, $w_{i,2}$ and~$w_{i,3}$ are both coloured red, while $w_{i,1}$ is coloured blue, a contradiction to the assumption that the assignment is NAE-satisfying.
    Thus, $c_i$ has at most $d$ neighbours of the other colour.
    Hence, the colouring is indeed a red-blue $d$-colouring of~$G$.
    This completes the proof.
\end{proof}

\noindent
We remark that for every $\Delta>2d+3$, {\sc $d$-Cut} is still NP-complete on claw-free graphs with maximum degree $\Delta$: after replacing each~$F_\ell$ as a copy of $H_{d,k_\ell,\Delta-1}$, instead of $H_{d,k_\ell,2d+2}$, the same hardness proof works.

\section{Conclusions}\label{s-con}

We proved both constant-time and \NP-completeness results for {\sc $d$-Cut} restricted to claw-free graphs of bounded maximum degree, leaving open only one case. We pose this open case as an open problem, replacing Open Problem~\ref{o-2}.

\begin{open}
    For every integer $d\geq2$, what is the computational complexity of {\sc $d$-Cut} on claw-free graphs of maximum degree~$\Delta=2d+2$?
\end{open}

\noindent
In our paper we also showed that our positive results for {\sc $d$-Cut} on claw-free graphs can be generalized to even $S_{1^t,\ell}$-free graphs.
A natural question, which we leave for future work, is whether our results can also be generalized to $S_{1,2,2}$-free graphs.
Here, $S_{1,2,2}$ is the graph obtained from the claw by subdividing two of its edges exactly once. 
In order to answer this question, we need some new arguments as our current proof techniques cannot be applied immediately. For instance,
in the proof of Theorem~\ref{t-star-free}, a key observation is that for each $u\in U$, $G[N_u]$ has no independent set of size~$t$, as otherwise we can find an induced $S_{1^t,\ell}$ by combining such an independent set with the path of length~$\ell$ between~$u$ and~$v$.
This observation no longer holds for the $S_{1,2,2}$ case.

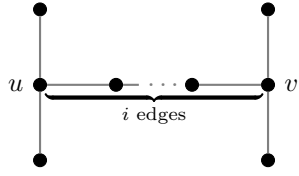
\begin{figure}[t]
    \centering
    \scalebox{1}{
\begin{tikzpicture}

\begin{scope}
	\node[vertex] (v1) at (0,2){};
	\node[vertex, label=left:$u$] (v2) at (0,1){};
	\node[vertex] (v3) at (0,0){};
	\node[vertex] (v4) at (3,2){};
	\node[vertex, label=right:$v$] (v5) at (3,1){};
	\node[vertex] (v6) at (3,0){};
	\node[vertex] (u1) at (1,1){};
	\node[vertex] (u2) at (2,1){};
	\draw[edge](v1)--(v2);
	\draw[edge](v2)--(v3);
	\draw[edge](v4)--(v5);
	\draw[edge](v5)--(v6);
	\draw[edge](v2)--(u1);
	\draw[edge](u2)--(v5);
	\draw[edge](u1)--(1.3,1);
	\node[color = gray](dots) at (1.65,1){$\dots$};	
	\draw[brace](2.925,0.875)--(0.075,0.875);
	\scriptsize
	\node[](i) at (1.5, 0.6){$i$ edges};
\end{scope}
\end{tikzpicture}
}
    \caption{The graph~$H_i^*$. Figure taken from~\cite{LMPS24}.}
    \label{f-h-graph}
\end{figure}

We recall that in particular, our results imply that {\sc $2$-Cut}, restricted to claw-free graphs of maximum degree~$p$, is constant-time solvable if $p\leq 5$ and \NP-complete if $p\geq 7$.
The latter result resolved an open case in the complexity classification of {\sc $d$-Cut} on $H$-free graphs, as the complexity of {\sc $2$-Cut} for claw-free graphs was not previously known. This brings us to the state-of-the-art summary for {\sc $d$-Cut} on $H$-free graphs, which we discuss and update below.

We first give some additional terminology.
For $r\geq 1$, we denote by~$P_r$ and~$C_r$, respectively, a path and a cycle on~$r$ vertices.
In addition, for $s\geq1$, we denote by $sP_r$ the disjoint union of $s$ copies of~$P_r$.
We denote by $S_{1,1,2}$ the graph obtained from a claw by subdividing one edge once, so that it has four edges.
Let~$H_1^*$ be a graph that looks like the letter ``H''.
It is obtained by taking the graph $2P_3$ and making the middle vertices, say $u$ and $v$, of each~$P_3$ adjacent to each other.
We obtain the graph $H_i^*$, for $i\geq 2$, by subdividing~$uv$ exactly ${i-1}$ times, see Figure~\ref{f-h-graph}.
The {\it girth} of a graph that is not a forest is the length of a shortest induced cycle in it.

We now present the updated state-of-the-art summary.
For every $d\geq 1$, {\sc $d$-Cut} is \NP-complete for graphs of girth at least $g$ for every fixed $g\geq 3$~\cite{FLPR23}, and thus for $C_r$-free graphs for every $r\geq 3$.
It is also known that {\sc $1$-Cut} is \NP-complete for $(H_1^*,\ldots,H_i^*)$-free  graphs for every $i\geq 1$~\cite{FLPR23}, $K_{1,4}$-free graphs~\cite{Ch84}, 
$(3P_5,P_{15})$-free graphs~\cite{LPR23} and $(3P_6,2P_7,P_{14})$-free graphs~\cite{LL23}. On the positive side, $1$-{\sc Cut} is polynomial-time solvable for $(sP_3+S_{1,1,2})$-free graphs~\cite{LMO25}, $(sP_3+P_4+P_6)$-free graphs~\cite{LMO25} and $(sP_3+P_7)$-free graphs~\cite{LMO25} for every $s\geq 1$.

We summarize the above results and our new result for $2$-{\sc Cut} on claw-free graphs in the following theorem, in which all unreferenced results for {\sc $d$-Cut} $(d\geq 2)$ are shown in~\cite{LMPS24}. 
For two graphs~$H$ and~$H'$, we write $H\ssi H'$ if~$H$ is an induced subgraph of~$H'$, and denote by $H+H'$ the disjoint union of a copy of~$H$ and a copy of~$H'$. 

\begin{theorem}\label{t-summary}
    Let $d\geq 1$ and let~$H$ be a graph.
    \begin{itemize}
        \item If $d=1$, then {\sc $d$-Cut} ({\sc Matching Cut}) on $H$-free graphs is 
        \begin{itemize}
            \item polynomial-time solvable if $H\ssi sP_3+S_{1,1,2}$, $sP_3+P_4+P_6$, or $sP_3+P_7$ for some $s\geq 1$;
            \item \NP-complete if $H\si K_{1,4}$, $P_{14}$, $2P_7$, $3P_5$, $C_r$ for some $r\geq 3$, or~$H_i^*$ for some $i\geq 1$. 
        \end{itemize}
        \item If $d\geq 2$, then {\sc $d$-Cut} on $H$-free graphs is
        \begin{itemize}
            \item polynomial-time solvable if $H\ssi  sP_1+P_3+P_4$ or $sP_1+P_5$ for some $s\geq 1$;
            \item \NP-complete if $H\si K_{1,3}$, $3P_2$, $C_r$ for some $r\geq 3$.
        \end{itemize}
    \end{itemize}
\end{theorem}

\noindent
It is known that for $d\geq 2$, {\sc $d$-Cut} is polynomial-time solvable for $(H+P_1)$-free graphs whenever {\sc $d$-Cut} is polynomial-time solvable for $H$-free graphs. This means that the cases $\{H+sP_1\; |\; s\geq 0\}$ are all polynomially equivalent. Hence, for $d\geq 2$, there are, due to our new result, now only three non-equivalent graphs $H$ left for which we do not know the computational complexity of {\sc $d$-Cut} on $H$-free graphs, namely $H=\{2P_4,P_6,P_7\}$.

\medskip
\noindent
\textit{Acknowledgements}.
The authors thank Ali Momeni and Hyunwoo Lee for providing some helpful insights about the proofs of Lemma~\ref{lem:cubic} and Theorem~\ref{t-star-free}, respectively.

\end{document}